\documentstyle[12pt]{article}

\newcommand{\nwc}{\newcommand}
\newcommand{\bes}{\begin{displaymath}}
\newcommand{\ees}{\end{displaymath}}
\newcommand{\be}{\begin{equation}}
\newcommand{\ee}{\end{equation}}

\newcommand{\bl}{{\bf l}}
\newcommand{\ba}{\begin{eqnarray}}
\newcommand{\ea}{\end{eqnarray}}
\newcommand{\bas}{\begin{eqnarray*}}
\newcommand{\eas}{\end{eqnarray*}}
\newcommand{\half}{\frac{1}{2}}
\newcommand{\al}{\alpha}
\newcommand{\bE}{{\bf E}}

\newcommand{\bI}{{\bf I}}

\newcommand{\bD}{{\bf D}}

\newcommand{\bt}{\beta}

\newcommand{\bbm}{{\bf m}}

\newcommand{\si}{\sigma}

\newcommand{\bB}{{\bf B}}

\newcommand{\vV}{{\bf V}}
\newcommand{\bZ}{{\bf Z}}

\newcommand{\ep}{\varepsilon}

\newcommand{\la}{\lambda}
\newcommand{\cL}{{\cal L}}

\newcommand{\Om}{\Omega}

\newcommand{\ga}{\gamma}
\newcommand{\bR}{{\bf R}}

\newcommand{\bk}{{\bf k}}

\newcommand{\bze}{{\bf 0}}

\newcommand{\kropa}{$_{\Box}$}
\newcommand{\bU}{{\bf U}}

\newcommand{\da}{\downarrow0}
\newtheorem{theorem}{Theorem}

\newtheorem{lemma}{Lemma}

\nwc{\m}{\mbox}
\nwc{\bm}{\boldmath}
\nwc{\ubm}{\unboldmath}
\nwc{\bmu}{\m{\bm $u$\ubm}}
\nwc{\bmx}{\m{\bm $x$\ubm}}
\nwc{\bx}{\bmx}
\newcommand{\by}{\m{\bm $y$\ubm}}
\nwc{\bmy}{\by}
\nwc{\bmv}{\m{\bm $v$\ubm}}
\newcommand{\vv}{\bmv}
\nwc{\cE}{{\cal E}}

\nwc{\beq}{\begin{eqnarray}}
\nwc{\eeq}{\end{eqnarray}}
\begin{document}
\title{Limit Theorems for  Motions
 in a Flow with a
Nonzero Drift.}
\author{Albert Fannjiang \\ Department of
Mathematics, University of California, Davis\and Tomasz Komorowski
\\ Institute of Mathemathics, Maria Curie-Sk\l odowska
University, Lublin}

\maketitle

{\bf Abstract} We establish  diffusion and fractional Brownian motion
approximations for motions in a
 Markovian Gaussian random field
with a nonzero mean.



 {\bf Abbreviated title} Limit theorems for motions in a random field.

\newpage

\section{Introduction.}
The simplest model of a motion of a tracer particle in a turbulent medium, e.g.  fluid, is given
by the stochastic differential equation
\be
\label{a1.0}
d\bmx(t)=\left(\vv+\ep\vV(t,\bmx)\right)dt+\sqrt{2\kappa}d\bB(t),
\ee
where the particle is located at $\bze$ when $t=0$.  $\vv\in R^d$ is a constant
vector representing the mean velocity of the fluid,
$\vV(t,\bmx)=(V_1(t,\bmx),\cdots,V_d(t,\bmx))$ describes the {\em velocity
fluctuations} and  is assumed to be a zero mean
 stationary random field.

$\bB(t)$ is a standard Brownian motion, whose presence accounts for the
existence of the molecular diffusivity in the medium given by $\kappa\geq 0$.
The parameter $\ep>0$ is supposed to indicate smallness of the fluctuation
amplitude in comparison with the mean drift. Turbulent  flows of that nature
have been extensively studied in Physics of Fluids   starting with the work of
the British physicist G.  Taylor in the 1920's (see \cite{T} and
 also \cite{frisch}
for more references on the subject).

The question is to determine the long time behavior of
the displacement $\bmx(t)$
 from the
statistics of the velocity and the size of the molecular diffusivity $\kappa$.
The effect of molecular diffusivity is in many cases  negligible if $\kappa$ is
small. We  therefore set $\kappa=0$ for simplicity of presentation. We
 assume
also that the mean velocity $\vv\not=0$. The case of vanishing mean velocity
 is presented by the authors in a separate paper (see \cite{fan-kom1,fan-kom}).

More specifically we set
\be
\label{971}
\by_\ep(t):=\bx\left(\frac{t}{\ep^2}\right)-\vv \frac{t}{\ep^2}.
\ee
We  then have
\[
\frac{d}{dt}\by_\ep(t)=\frac{1}{\ep}\vV
\left(\frac{t}{\ep^2},\vv \frac{t}{\ep^2}+\by_\ep(t)\right).
\]

 A theorem describing the limiting behavior of the particle,
proven in various forms in \cite{khasminskii,kp,kom}, says that $\by_\ep(t)$,
$t\geq0$, considered as continuous trajectory stochastic processes tend weakly,
as $\ep\da$, to a Brownian Motion whose covariance matrix is given by the
so-called {\em Taylor - Kubo formula}
\be
\label{TK}
D_{ij}^*=\int\limits^\infty_0\left\{\bE[V_i(t,\vv t)V_j(0,\bze)]+
\bE[V_i(t,\vv t)V_j(0,\bze)]\right\}dt
\ee
(cf. \cite{K}).

The  version of the above theorem proven in \cite{khasminskii} for $\vv=0$,
requires, among others, the  strong mixing assumption for the velocity field in
the temporal variable. The key point of the argument is an  approximation of the
forward particle location valid for sufficiently large times that only takes
into account  the information available prior to a given instant, the so-called
{\em frozen path approximation}. The temporal mixing property of the field
guarantees then that
 the future particle displacement is almost independent of the past information thus,
 in consequence the diffusion approximation holds. This type of theorem can be also
shown for a field with $\vv\not=\bze$ provided that its fluctuation is strongly
mixing in space (see \cite{kp}). In this case one can project the spatial domain
into time domain via the change of variables $\bx:=\bx-\vv t$, cf. (\ref{971}),
and then, using the fact that the mean drift size dominates, for small $\ep$,
the amplitude of fluctuations it is possible to translate the spatial mixing
property into temporal one for the velocity field considered in the new
coordinates and, as a result, the argument of \cite{khasminskii} can be applied.

 In this paper
we set out to give a description of the limiting behavior of scaled trajectories
for  a class of Gaussian fluctuation fields with {\em power-law spectrum}. This
family of  fields appears frequently in the mathematical theory of turbulence
and is motivated by the seminal work of Kolmogorov \cite{kol}, see also
\cite{am} for the shear layer flow case. In many interesting cases the
fluctuation fields with that kind of spectrum
 do not satisfy the mixing condition either in time or space, thus the
existing versions of limit theorems  are not applicable.

In what follows we  assume that $\vV(t,\bmx)$, $(t,\bmx)\in R\times R^d$, is a
$d$-dimensional, time-stationary, space-homogeneous, time-Markovian, isotropic
Gaussian velocity field  defined over a probability space
 $(\Om,{\cal V},P)$
whose two-point covariance matrix $\bR=[R_{ij}]$ is given by the Fourier
transform:
\be
\label{a1.3}
R_{ij}(t,\bmx)=\bE\left[V_i(t,\bmx)V_j(0,\bze)\right]=
\int\limits_{R^d}c_0(\bk\cdot\bx)e^{-|\bk|^{2\bt}t}
\hat{R}_{i,j}(\bk)d~\bk.
\ee
Here and in the sequel we let
 $c_0(\phi):=\cos\phi,~~c_1(\phi):=\sin\phi$. $\bE$
stands for the mathematical expectation calculated with respect to
 $P$. The
spatial spectral density of the field $\hat{\bR}=[\hat{R}_{i,j}]$ is given by
\be
\label{060103}
\hat{R}_{i,j}(\bk)=\frac{a(|\bk|)}{|\bk|^{2\al+d-2}}\left(\delta_{i,j}-\frac{k_ik_j}{|\bk|^2}\right),
\ee
where $a:[0,+\infty)\rightarrow R_+$, the so-called {\em ultraviolet
 cut-off}, is
a certain compactly supported continuous  function. The factor
$\delta_{i,j}-k_ik_j/|\bk|^2$ in (\ref{060103}) ensures that the velocity field
is incompressible, i.e.
$\nabla\cdot\vV(t,\bx):=\sum\partial_iV_i(t,\bx)\equiv0$. The function
$\exp(-|\bk|^{2\beta}t)$ in (\ref{a1.3}) is called the {\em time correlation
function} of the velocity $\vV$. The  spectral density $\hat{\bR}(\bk)$ is
integrable over $\bk$
 for $\alpha<1$. The ultraviolet
cut-off  is then needed  to avoid divergence of the integral over $|\bk|$. The
parameter $\alpha$ is directly related to the decay exponent of
 $\bR$. Namely,
$\bR(\bmx)\sim |\bmx|^{\al-1}$ for $|\bmx|\gg1$. As $\alpha$ increases to one,
the spatial decay exponent of $\bR$
 decreases to zero and in consequence the strength of spatial
correlation of modes increases.  For $\alpha\geq 1$, the  spectral
 density
ceases to be integrable and needs an infrared cut-off, i.e. the origin should
lie outside of the support of $a$. On the other hand, in order to apply the
results of \cite{kp} or \cite{khasminskii} one needs to assume that $\bt\leq0$.
Otherwise, i.e. when $\beta>0$, the velocity field lacks the spectral gap and
thus the strong mixing property (cf. \cite{R}). In light of the foregoing
discussion we therefore restrict our attention to the case $\alpha<1, ~\beta
>0$, which corresponds to the velocity field having arbitrarily long scales but not
 the strong mixing property.

In addition, the particular form of the spectrum of the field allows us to
assume without any loss of generality
 that it is jointly continuous in both $t$  and $\bmx$ and is of
$C^\infty$ class in $\bmx$ a.s. (see e.g. \cite{adler}).

One can then pose  the following question: what is the region in the $(\alpha,
\beta)$ plane where the classical turbulent diffusion theorem, with the Taylor - Kubo formula (\ref{TK}),
holds? It is easy to find the necessary condition by imposing the convergence of
the integral appearing in (\ref{TK}). We have
\be
\label{0220}
D^*_{ij}=\int\limits^\infty_0 R_{ij}(t,\vv t)~dt
=\int\limits_{R^d}\left(\delta_{ij}-{k_ik_j\over |\bk|^2}\right)
\frac{a(|\bk|)}{|\bk|^{2\al+d-2}}\int\limits^\infty_0
\exp(-|\bk|^{2\beta}t)c_0(k_{\vv} t)~dt~ d\bk,
\ee
where $k_{\vv}:=\bk\cdot\vv$. After a straightforward calculation of the
integral over $t$  one deduces that the  rightmost part of (\ref{0220}) equals
\be
\label{a1.10}
\int\limits_{R^d}{a(|\bk|)|\bk|^{2\bt}\over
|\bk|^{2\al+d-2}(|\bk|^{4\bt}+k_{\vv}^2)}
\left(\delta_{ij}-{k_ik_j\over |\bk|^2}\right)
~d\bk.
\ee
Elementary calculations (cf. the computations made
 after (\ref{281101}) and (\ref{40'}) below)
show that the finiteness of the expression (\ref{a1.10}) leads to the conditions
\be
\label{a1.5}
\alpha+\beta<1~\mbox{ when }\bt<1/2
\ee
or
\be
\label{a11.5}
\alpha<1/2\leq\beta<1.
\ee

It turns out that these conditions are also sufficient to claim that the
particle trajectory approximates the Brownian motion in the sense of weak
convergence of continuous path stochastic processes. This   is essentially the
contents of our Theorem \ref{theorem1}  below.

One can also ask what happens in the case when conditions
 (\ref{a1.5}), (\ref{a11.5}) are  violated. The fact that
 the integral in (\ref{0220}) diverges suggests  ``slower''
than $\ep^{-2}$ scaling of the temporal variable  in   (\ref{971}). We define
\be
\label{120801}
\by_\ep(t):=\bx\left(\frac{t}{\ep^{2\delta}}\right)-
\vv\frac{t}{\ep^{2\delta}}
\ee
for a certain suitably chosen $\delta\in(1/2,1)$. The limit of the stochastic
processes given by the trajectories turns out then to be  a fractional Brownian
motion. This result is contained in Theorem~\ref{theorem2}.

It is interesting, in our judgement, to compare the results of the present
article with those of \cite{fan-kom1} and \cite{fan-kom}. It was proven there
that the diffusive regime for the motion in the field without the drift, i.e.
$\vv=0$, comprises of $(\al,\bt)$ satisfying $\al+\bt<1$ (see \cite{fan-kom}).
When, on the other hand, $\al+\bt>1$, the properly scaled trajectories
approximate a fractional Brownian motion (see \cite{fan-kom1}). As a result we
see that inthe present setting the diffusivity regime is enlarged for $\bt<1/2$.
We call this phenomenon, after \cite{am1}, the {\em sweeping effect}. It can be
explained intuitively as follows. For a given wavenumber of magnitude $|\bk|$
one can see from (\ref{a1.3}) that the {\em relaxation time} needed for a
significant temporal decorrelation of the fluctuation field is of order
$\tau\sim|\bk|^{-2\bt}$. On the other hand the {\em sweeping time} necessary for
a significant decorrelation related to the mean drift convection is of order of
magnitude $\tau_{\vv}\sim|\bk|^{-1}$. It is much shorter than the relaxation
time when $\bt\geq1/2$. In that regime therefore the mean drift is responsible
for the chaotic behavior of the velocity field and in consequence the diffusive
approximation holds when $\al\leq 1/2$. On the other hand for $\bt<1/2$ the
effect of the mean drift is negligible and the limiting particle motion is the
result of a subtle balance between the spatial and temporal chaotic properties
of the field as described in \cite{fan-kom1,fan-kom}.

To simplify the presentation we focus our attention only on the case when
$\al,\bt$ satisfy in addition
\be
\label{120802}
\al+2\bt<2.
\ee
This condition can be proven inessential for the validity of our results. The
more technically involved case when $\al+2\bt\geq2$
  will be the subject of a subsequent article.

We summarize our results in the form of the following two theorems.
\begin{theorem}
\label{theorem1}
Suppose that
\begin{itemize}
\item[1)] $\vv\not=0$ and {\em (\ref{120802})} holds.
\item[2)] $\vV(t,\bx)$, $(t,\bx)\in R\times R^d$, is
a  zero mean, Gaussian velocity field whose correlation is given by {\em
(\ref{a1.3})} and {\em (\ref{060103})}.
\end{itemize}

Then,
 for $(\al,\bt)$ as specified by conditions {\em (\ref{a1.5})} and {\em (\ref{a11.5})}  the scaled
  trajectories given by (\ref{971}) converge weakly, when $\ep\da$, as continuous
  trajectory stochastic processes
 to a Brownian motion whose
covariance matrix $\bD_0$ is given by  {\em (\ref{TK})}.
\end{theorem}

\begin{theorem}
\label{theorem2}
Assume  conditions {\em 1)} and {\em 2)} of Theorem {\em \ref{theorem1}}. Then
for $1>\al$, $\bt>0$ such that
\be
\label{a111.5}
\alpha+\beta\geq 1~~\mbox{when }1/2\leq\al<1
\ee
one can find a unique $\delta_{\al,\bt}\in(1/2,1)$ and
$H=\frac{1}{2\delta_{\al,\bt}}$ such that the scaled trajectories given by {\em
(\ref{120801})}  converge weakly to a fractional Brownian Motion $\bB_H(t)$,
$t\geq0$, i.e. the unique $H$-self-similar, Gaussian process with stationary
increments (see  {\em \cite{st}}). We have
\be
\label{010514}
\bE\left[\bB_H(t)\otimes \bB_H(t)\right]=\bD_{\al,\bt} t^{2H},
\ee
with
\be
\label{010513}
\bD_{\al,\bt}=\int\limits_{R^d}\Gamma_{\al,\bt}
\left(\bI-\frac{\bk\otimes\bk}{|\bk|^2}\right)
\frac{a(0)d\bk}{|\bk|^{d-1}}.
\ee
where $\Gamma_{\al,\bt}$ equals
\[
\frac{\mbox{e}^{-|\bk|^{2\bt}}-1+|\bk|^{2\bt}}{|\bk|^{2\al+4\bt-1}}
\qquad\mbox{when  }\bt<1/2,
\]
\[
\frac{1-c_0(k_{\vv})}{k_{\vv}^2|\bk|^{2\al-1}}
\qquad\mbox{when  }\bt>\half,
\]
\[
\frac{(|\bk|^2-k_{\vv}^2)\left(\mbox{e}^{-|\bk|}c_0(k_{\vv})-1+
|\bk|\right)-2|\bk|k_{\vv}\left(\mbox{e}^{-|\bk|} c_1(k_{\vv})-k_{\vv}\right])}
{(|\bk|^2+k_{\vv}^2)|\bk|^{2\al-1}}
\qquad\mbox{when  }\bt=\half.
\]
In addition, $\delta_{\al,\bt}$ equals
\be
\label{01051}
\frac{\bt}{\al+2\bt-1}\qquad \mbox{when }\bt<1/2,
\ee
\be
\label{01052}
\frac{1}{2\al}\qquad\mbox{when }\bt\geq1/2.
\ee

\end{theorem}


\section{Preliminaries.}The Spectral Theorem
for real vector valued random fields (cf. e.g. \cite{adler}) implies that
 there exist two
independent, identically distributed, real vector valued Gaussian spectral
measures $\hat{\vV}_{l}(t,\cdot)=
(\hat{V}_{l,1}(t,\cdot),\cdots,\hat{V}_{l,d}(t,\cdot))$, $l=0,1$, such that
\be
\label{981}
\vV(t,\bx)=\int \hat{\vV}_0(t,\bx,d\bk),
\ee
where
\[
\hat{\vV}_0(t,\bx,d\bk):=c_0(\bk\cdot\bx)\hat{\vV}_0(t,d\bk)+
c_1(\bk\cdot\bx)\hat{\vV}_1(t,d\bk).
\]
In the sequel (cf. e.g. Lemma \ref{lemma2}) we  also deal with the spectral
measure given by
\[
\hat{\vV}_1(t,\bx,d\bk):=-c_1(\bk\cdot\bx)\hat{\vV}_0(t,d\bk)+
c_0(\bk\cdot\bx)\hat{\vV}_1(t,d\bk).
\]

The field $\vV$ is {\em Markovian}
 in the following sense. For any function
$\psi\in{\cal S}(R^d,R)$ and $i\in\{1,\cdots,d\}$  we have
\be
\label{060110}
\label{1.10}
\bE\left[\int\psi(\bk)\hat{V}_{l,i}(t,\bx,d\bk)\left|\right.{\cal
    V}_{-\infty,s}\right]= \int
\mbox{e}^{-|\bk|^{2\bt}(t-s)}\psi(\bk)\hat{V}_{l,i}(s,\bx,d\bk),~~l=0,1.
\ee
Here,  for any $a\leq b$, ${\cal
    V}_{a,b}$ denotes the $\si$-algebra generated by $\vV(t,\bx)$ with
$(t,\bx)\in[a,b]\times R^d$ and $\hat{V}_{l,i}(s,\bx,d\bk)$ is
the relevant component of $\hat{\vV}_l(s,\bx,d\bk)$. In what follows
 we
also write $L^2_{a,b}$ for $L^2(\Om,{\cal
    V}_{a,b},P)$.

We now introduce the key concepts appearing in the proofs of Theorems
\ref{theorem1} and \ref{theorem2}, namely $\chi_{\la}$, the {\em
$\la$-corrector} and $\bU_\la$, the {\em $\la$-convector}. We set
\be
\label{123}
\chi_{\la}(t,\by):=\int\limits_t^{+\infty}\mbox{e}^{-\lambda(s-t)}
\bE\left[\vV(s,\vv s+\by)\left|\right.{\cal
  V}_{-\infty,t}\right]ds,
\ee
\be
\label{122}
\bU_{\la}(t,\by):=\vV(t,\vv t+\by)\cdot\nabla\chi_{\la}(t,\by).
\ee
The differentiation used in (\ref{122}) is understood in the mean
square sense.

The procedure defining $\chi_{\la}$ and $\bU_{\la}$ requires some justification
at least to guarantee smoothness of the corrector field appearing on the right
hand side of (\ref{122}). This technical point is explained by the following.
\begin{lemma}
\label{lemma2}
For arbitrary $\la>0$ the fields
 $\chi_{\la}$, $\bU_{\la}$ defined
above are  $L^p$ integrable for any $p\geq 1$. $L^p$ norms of those fields have
the following hypercontractivity property:
\be
\label{061201}
\left(\bE|\chi_{\la}(t,\by)|^p\right)^{1/p}\leq
C\left(\bE|\chi_{\la}(t,\by)|^2\right)^{1/2},
\ee
\be
\label{061201a}
\left(\bE|\bU_{\la}(t,\by)|^p\right)^{1/p}\leq
C\left(\bE|\bU_{\la}(t,\by)|^2\right)^{1/2},
\ee
where the constant $C>0$ does not depend on $\la$.

One can also
 select  modifications of $\chi_{\la}$ and $\bU_{\la}$
which are jointly   stationary in the strict sense, continuous in both $(t,\bx)$
and $C^\infty$  in $\bx$.

The spectral representation of $\chi_\la$ is given by
\[
\chi_\la(t,\by)=\chi_\la^{(0)}(t,\by)+\chi_\la^{(1)}(t,\by)
\]
with
\be
\label{4161}
\chi_\la^{(i)}(t,\by):=
\int C_i(\bk,\la)\hat{\vV}_i(t,\vv t+y,d\bk)~i=1,2
\ee
when
\be
\label{4162}
C_1(\bk,\la):=\frac{|\bk|^{2\bt}+\la}{(|\bk|^{2\bt}+\la)^2+ k_{\vv}^2}
\ee
and
\be
\label{41610}
C_2(\bk,\la):=
\frac{k_{\vv}}{(|\bk|^{2\bt}+\la)^2+k_{\vv}^2}
\ee
\end{lemma}

The proof of the lemma is fairly standard. We sketch here only its main points
referring the reader interested in details to the relevant literature.
  $L^p$ integrability and inequalities (\ref{061201}), (\ref{061201a}) for the
fields in question  follow from   the hypercontractivity property of $L^p$ norms
in Gaussian measure spaces related to the velocity field $\vV$ (see
 e.g. Theorem
5.1 and its corollaries  in \cite{janson}).

The existence of  regular versions of the fields follows
from the fact that there exists $h>0$ such that for any integer $N>0$ we can find a constant $C_{N,h}>0$ depending on $h,N$ only
for which
\[
\sum\limits_{|\bbm|=N}\bE|D^\bbm\vV(t,\bx)-D^\bbm\vV(s,\by)|^2\leq
C_{N,h}\left[|t-s|^h+|\bx-\by|^2\right]
\]
for all $(t,\bx),(s,\by)\in R\times R^d$. Here for any integral
 multiindex
$\bbm=(m_1,\cdots,m_d)$ we define
$D^{\bbm}:=\partial_{x_1}^{m_1}\cdots\partial_{x_d}^{m_d}$. We can find the
modifications of the fields with required degree of regularity using the
aforementioned Theorem 5.1 of ibid in conjunction with Kolmogorov's classical
result on the existence of continuous trajectory modification of a stochastic
process (the random field version of that criterion can be found in e.g.
\cite{adler} as the Corollary to Theorem 3.2.5).

Joint stationarity of the fields in question can easily be verified by an
application of the results of \cite{port-stone}.

\section{Proof of Theorem \ref{theorem1}.}
We define a {\em scaled corrector along  path} by
\be
\label{126}
\chi_{\ep}(t):=\ep\chi_{\ep^2}(\frac{t}{\ep^{2}},\by_\ep(t))
\ee
and a {\em  scaled convector along  path} by
\be
\label{129}
\bU_{\ep}(t):=\bU_{\ep^2}(\frac{t}{\ep^{2}},\by_\ep(t)).
\ee

Thanks to the divergence free structure of the velocity field $\vV$ the
processes $\chi_{\ep}(t)$, $t\geq0$, and $\bU_{\ep}(t)$, $t\geq0$, are jointly
strictly stationary  (see
 Theorem
2, p. 500, of \cite{port-stone}).

The importance of the concepts of a corrector
and convector is highlighted by
the following lemma.
\begin{lemma}
\label{lemma1}
We have
\be
\label{124}
\by_\ep(t)=
-\chi_{\ep}(t)+\chi_{\ep}(0)+\int\limits_0^t\chi_{\ep}(s)~ds+
\int\limits_0^t\bU_{\ep}(s)~ds+
M_\ep(t).
\ee
Here $M_\ep(t)$, $t\geq0$, is a Brownian motion whose covariance matrix is given
by $\bD_\ep $ where
\be
\label{061104}
\bD_\ep=\int\limits_{R^d}\frac{|\bk|^{2\bt}a(|\bk|)}
{|\bk|^{2\al+d-2}\left[(|\bk|^{2\bt}+\ep^2)^2+k_{\vv}^2\right]}\left(\bI-\frac{\bk\otimes\bk}{|\bk|^2}\right)d\bk.
\ee
In addition, there exists $\ga>0$ such that
\be
\label{08121}
\lim\limits_{\ep\da} \ep^{-\ga}\bE|\chi_{\ep}(t)|^2=0 \qquad\mbox{and}\qquad
\label{08122}
\lim\limits_{\ep\da}  \ep^{-\ga}\bE|\bU_{\ep}(t)|^2=0.
\ee
\end{lemma}

{\bf Proof.} First we show (\ref{124}). We start with  the calculation of
 the so-called {\em pseudogenerator} of the path
corrector, i.e. the process given by (cf. Chapter 3 of \cite{kushner})
\be
\label{010516}
L\chi_\ep(t):=
\lim\limits_{\delta\da}\frac{\bE\left[\chi_\ep(t+\delta)
\left|\right.{\cal
  V}_{-\infty,t/\ep^2}\right]-\chi_\ep(t)}{\delta},
\ee
where the limit is understood in the $L^1$ sense.

It is easy to observe  via Taylor's expansion of
$\chi_\ep(\frac{t+\delta}{\ep^2},\by)$ about $\by_\ep(t)$ that
\be
\label{9112}
\chi_\ep(t+\delta)=
\ep\chi_{\ep^2}\left(\frac{t+\delta}{\ep^2},\by_\ep(t)\right)+
R_\delta(t),
\ee
where
\[
R_\delta(t)=
\bU_{\ep^2}\left(\frac{t}{\ep^2},\by_{\ep}(t)\right)\delta+o(\delta).
\]
Here $o(\delta)$ is understood as a function with takes values in the space
random vectors with $L^p$ integrable components for an arbitrarily chosen $p\geq
1$ and such that $o(\delta)/\delta$ vanishes in the $L^p$ sense as $\delta\da$.

Conditioning the first term on the right hand side with respect to the
$\si$-algebra ${\cal
  V}_{-\infty,t/\ep^2}$ we deduce from (\ref{123}) and
(\ref{122}) that
\be
\label{995}
\ep\bE\left[\chi_{\ep^2}\left(\frac{t+\delta}{\ep^2},\by_\ep(t)\right)\left|\right.{\cal
  V}_{-\infty,t/\ep^2}\right]=
\ep\int\limits_{t+\delta/\ep^2}^{+\infty}
\mbox{e}^{-\la(s-(t+\delta)/\ep^2)}\bE\left[\vV(s,\vv s+\by_\ep(t))
\left|\right.{\cal
  V}_{-\infty,t/\ep^2}\right].
\ee
The last expression, after  applying Taylor's expansion this time about
$t/\ep^2$, can be rewritten as
\be
\label{994}
\chi_\ep(t)+\delta\left[\chi_\ep(t)-\frac{1}{\ep}\vV\left(\frac{t}{\ep^2},\vv\frac{t}{\ep^2}+
\by_\ep(t)\right)\right]+o(\delta).
\ee

In consequence,
\be
\label{41510}
L\chi_\ep(t)= \chi_\ep(t)-
\frac{1}{\ep}\vV\left(\frac{t}{\ep^2},\vv\frac{t}{\ep^2}+
\by_\ep(t)\right)+
\bU_\ep(t).
\ee
Using  the results of Chapter 3 of \cite{kushner}
 we obtain that
\be
\label{3111}
M_\ep(t)=\chi_\ep(t)-\int\limits_0^t L\chi_\ep(s)~ds
\ee
is a continuous trajectory vector valued martingale. We write
$M_{\ep}=(M_{\ep,1},\cdots,M_{\ep,d})$; a similar convention will be used for
other vector valued fields appearing with subscript $\ep$.

We notice that thanks to the aforementioned results of \cite{kushner} the joint
quadratic variation of martingales $M_{\ep,i}(t),~M_{\ep,j}(t)$, $t\geq0$,
$i,j=1,\cdots,d$, equals
\be
\label{061602}
\langle M_{\ep,i},M_{\ep,j}\rangle_t=\int\limits_0^tL\left\{M_{\ep,i}(s)M_{\ep,j}(s)\right\}ds.
\ee
After  elementary calculations we get
\be
\label{061601}
L\left\{M_{\ep,i}(t)M_{\ep,j}(t)\right\}=
\lim\limits_{\delta\da}\frac{ \bE\left[M_{\ep,i}(t+\delta)M_{\ep,j}
(t+\delta) \left|\right.{\cal
  V}_{-\infty,t/\ep^2}\right]-
M_{\ep,i}(t)M_{\ep,j}(t)}{\delta}
\ee
\[
=\lim\limits_{\delta\da}\frac{ \bE\left[\chi_{\ep,i}(t+\delta)\chi_{\ep,j}
(t+\delta) \left|\right.{\cal
  V}_{-\infty,t/\ep^2}\right]- \chi_{\ep,i}(t)\chi_{\ep,j}
(t)}{\delta}-\left[\chi_{\ep,i}(t) L\chi_{\ep,j}(t) +
\chi_{\ep,j}(t)L \chi_{\ep,i}(t)\right].
\]
Using  expansion (\ref{9112}) to represent $\chi_{\ep}(t+\delta)$ in
(\ref{061601}) we obtain that the utmost right hand side of (\ref{061601})
equals
\be
\label{4152}
L_{\ep,i,j}(t;\by_\ep(t))+U_{\ep,i}(t)\chi_{\ep,j}(t)
+\chi_{\ep,i}(t)U_{\ep,j}(t)-\chi_{\ep,i}(t) L\chi_{\ep,j}(t)
-\chi_{\ep,j}(t)L\chi_{\ep,i}(t),
\ee
where
\be
\label{993}
L_{\ep,i,j}(t;\by):=\ep^2\lim\limits_{\delta\da}
\frac{ \bE\left[\chi_{\ep,i}(\frac{t+\delta}{\ep^2},\by)
\chi_{\ep,j}
(\frac{t+\delta}{\ep^2},\by)\left|\right.  {\cal
  V}_{-\infty,t/\ep^2}\right]-
\chi_{\ep,i}(\frac{t}{\ep^2},\by)
\chi_{\ep,j}
(\frac{t}{\ep^2},\by)}{\delta}.
\ee
The corrector
 fields appearing in (\ref{993}) are zero mean  Gaussian therefore
\[
\bE\left[\chi_{\ep,i}\left(\frac{t+\delta}{\ep^2},\by\right)
\chi_{\ep,j}
\left(\frac{t+\delta}{\ep^2},\by\right)\left|\right.  {\cal
  V}_{-\infty,t/\ep^2}\right]=m_{i,j}-\bE m_{i,j}+\bE\left[\chi_{\ep,i}(0)\chi_{\ep,j}(0)\right]
\]
where
\[
m_{i,j}=\ep^2\bE\left[\chi_{\ep,i}
(\frac{t+\delta}{\ep^2},\by) \left|\right.  {\cal
  V}_{-\infty,\frac{t}{\ep^2}}\right]
\bE\left[\chi_{\ep,j}\left(\frac{t+\delta}{\ep^2},\by\right)\left|\right.  {\cal
  V}_{-\infty,t/\ep^2}\right].
\]
With the help of (\ref{995}) and (\ref{994}) we deduce that
\be
\label{4151}
L_{\ep,i,j}(t;\by_\ep(t))=B_{\ep,i,j}(t)-\bE B_{\ep,i,j}(t)
\ee
where
\[
B_{\ep,i,j}(t):=
\left[\chi_{\ep,i}(t)-\frac{1}{\ep}
V_i\left(\frac{t}{\ep^2},\vv\frac{t}{\ep^2}+
\by_\ep(t)\right)\right]\chi_{\ep,j}(t)+
\left[\chi_{\ep,j}(t)-\frac{1}{\ep}
V_j\left(\frac{t}{\ep^2},\vv\frac{t}{\ep^2}+
\by_\ep(t)\right)\right]\chi_{\ep,i}(t).
\]
Combining (\ref{4151}) with (\ref{4152}) we deduce,
 using (\ref{41510}), that
\be
\label{4153}
L\left\{M_{\ep,i}(t)M_{\ep,j}(t)\right\}=-\bE B_{\ep,i,j}(t)=-\bE
B_{\ep,i,j}(0),
\ee
with the last equality following from stationarity of the relevant processes
along the particle trajectory that holds due to the results of
\cite{port-stone}. A direct computation, using the formula for the expectation
of a product of Gaussian random variables, shows that the right hand side of
(\ref{4153}) equals $D_{\ep,i,j}$, where $\bD_\ep:=[D_{\ep,i,j}]$ is given by
(\ref{061104}). This concludes the proof of (\ref{124}).

Now we prove  formula (\ref{08121}). Without any loss of generality we assume
that $\vv=(1,0,\cdots,0)$.

We then have
\be
\label{281101}
\bE\left|\chi_\ep(t)\right|^2=\bE\left|\chi_\ep(0)\right|^2\leq
C\ep^2\int\limits_{|\bk|\leq K}
\frac{d\bk}{\left[(\ep^2+|\bk|^{2\bt})^2+
k_1^2\right]|\bk|^{2\al+d-2}}.
\ee
If we represent $\bk=(k_1,\bl)$, where $\bl\in R^{d-1}$ the rightmost
 part of (\ref{281101}) can be rewritten as
\be
\label{010101}
C\ep^2\mathop{\int\int}\limits_{\sqrt{k_1^2+|\bl|^2}\leq K}
\frac{dk_1d\bl}{\left\{[\ep^2+(k_1^2+|\bl|^2)^{\bt}]^2+k_1^2\right\}
(k_1^2+|\bl|^2)^{\al+(d-2)/2}}
\ee
\[
\leq C'\mathop{\int\int}\limits_{\sqrt{k_1^2+l^2}\leq K}
\frac{\ep^2 l^{d-2}dk_1dl}
{\left\{[\ep^2+(k_1^2+l^2)^{\bt}]^2+k_1^2\right\}
(k_1^2+l^2)^{\al+(d-2)/2}}.
\]
After changing variables in the last integral according to the rule
$l=r\cos\phi$, $k_1=r\sin\phi$ we infer that the leftmost part of (\ref{281101})
can be estimated from above by
\[
C\int\limits_0^K\int\limits_0^{\pi/2}\frac{\ep^2drd\phi}
{[(\ep^2+r^{2\bt})^2+r^2\phi^2]r^{2\al-1}}.
\]
Here $C$ denotes a generic constant independent of $\ep$.

After an elementary calculation we find that the last integral
equals
\[
\int\limits_0^K\frac{\ep^2}{r^{2\al}(\ep^2+r^{2\bt})}
\arctan\left(\frac{\pi r}{2(\ep^2+r^{2\bt})}\right)dr.
\]
For $\al,\bt$ as specified by conditions (\ref{a1.5}), (\ref{a11.5}) the
Dominated Convergence Theorem implies that this expression vanishes as $\ep\da$.

To show (\ref{08122}) we first observe  that
\be
\label{221}
\bE\left|\bU_{\la}(0,\bze)\right|^2\leq
C\bE\left|\vV(0,\bze)\right|^2
\bE\left|\nabla\chi_{\la}(0,\bze)\right|^2.
\ee
The gradient of the corrector satisfies
\be
\label{40'}
\bE\left|\nabla\chi_{\la}(0,\bze)\right|^2\leq
\int\limits_{|\bk|\leq K}\frac{|\bk|^2}
{\left[(\lambda+|\bk|^{2\bt})^2+k_1^2\right] |\bk|^{2\al+d-2}}~d\bk.
\ee
Estimating precisely in the same way as we did for (\ref{281101}) we infer that
the right hand side of (\ref{40'}) is less than or equal to a constant times
\[
\int\limits_0^K\frac{1}{r^{2\al-2}(\la+r^{2\bt})}
\arctan\left(\frac{\pi r}{2(\la+r^{2\bt})}\right)dr.
\]
This expression
 remains bounded as $\la\da$ if only $\al+2\bt<2$.
Consequently,
\be
\label{222}
\limsup\limits_{\la\da}\bE\left|\bU_{\la}(0,\bze)\right|^2<+\infty.
\ee
(\ref{222}) holds also for any $p$-th absolute moment of $\bU_{\la}$ with $p>0$
thanks to Lemma \ref{lemma2}.

We also have
\[
\bE\left|\int\limits_0^t
\bU_{\ep^2}\left(\frac{s}{\ep^{2}},\bmy_\ep(s)\right)ds\right|^2=
2\int\limits_0^tds\int\limits_0^s\bE\left\{
\bU_{\ep^2}\left(\frac{s}{\ep^{2}},\bmy_\ep(s)\right)\cdot
\bU_{\ep^2}\left(\frac{s_1}{\ep^{2}},\bmy_\ep(s_1)\right)\right\}ds_1
\]
\be
\label{31.12}
=2\int\limits_0^tds\int\limits_0^s\bE\left\{\bU_{\ep^2}
\left(\frac{s}{\ep^{2}},\bmy_\ep(s_1)\right)\cdot
\bU_{\ep^2}\left(\frac{s_1}{\ep^{2}},\bmy_\ep(s_1)\right)\right\}ds_1
\ee
\[
+\frac{2}{\ep}\int\limits_0^tds\int\limits_0^sds_1
\int\limits_{s_1}^s\bE\left\{\left(\vV\left(\frac{s_2}{\ep^{2}},
\vv\frac{s_2}{\ep^2}+\bmy_\ep({s_2})\right)\cdot
\nabla\bU_{\ep^2}\left(\frac{s}{\ep^{2}},\by_\ep({s_2})\right)\right)\cdot
\bU_{\ep^2}\left(\frac{{s_1}}{\ep^{2}},\bmy_\ep(s_1)\right)\right\}d{s_2}
\]
The first term of (\ref{31.12}) equals
\be
\label{01021}
2\int\limits_0^tds\int\limits_0^s
\bE\left\{\bU_{\ep^2}\left(\frac{s_1}{\ep^{2}},\bze\right)
\cdot\bU_{\ep^2}(0,\bze)\right\}d{s_1}
\ee
due to the stationarity   of the path convector.

Thanks to the incompressibility of the velocity field we have $\bE\bU_{\la}
(t,\bmy)=\bze$. That and  well known properties of
 conditional expectations of  second degree polynomial-like functionals
of a Gaussian field  then imply that
\be
\label{01023}
\bU_{\la} (t,\bmy;0)=\vV(t,\vv t+\bmy;0)\cdot \nabla\chi_\la(t,\by;0).
\ee
Here $\vV(\cdot,\cdot;0)$, $\nabla\chi_\la(\cdot,\cdot;0)$ denote the orthogonal
projections of $\vV(\cdot,\cdot)$ and $\nabla\chi_\la(\cdot,\cdot)$ onto
$L^2_{-\infty,0}$.

In consequence, (\ref{01021}) equals
\be
\label{1.5}
2\int\limits_0^tds\int\limits_0^s
\bE\left\{\left(\vV\left(\frac{{s_1}}{\ep^{2}},
\vv\frac{s_1}{\ep^2};0\right)
\cdot\nabla\chi_{\ep^2}\left(\frac{s_1}{\ep^{2}},\bze;0\right)\right)
\cdot\bU_{\ep^2}(0,\bze)\right\}d{s_1}.
\ee
Using  the Cauchy Schwarz inequality we estimate (\ref{1.5}) by
\be
\label{0303}
C_1(\bE|\bU_{\ep^2}(0,\bze)|^4)^{1/4}(\bE|\vV(0,\bze)|^4)^{1/4}
\left(\int\limits_0^tds\int\limits_0^s
\bE\left|\nabla\chi_{\ep^2}
\left(\frac{{s_1}}{\ep^{2}},\bze;0\right)\right|^2d{s_1}
\right)^{1/2}.
\ee

The fourth moment of $\bU_{\ep^2}$ remains bounded as $\ep\da$, see (\ref{222}).
Therefore, the expression in (\ref{0303}) can be estimated  by a constant times
the last factor of (\ref{0303}). The latter, with the help of (\ref{1.10}) and
(\ref{40'}), equals
\be
\left[\int\limits_0^t\int\limits_{|\bk|\leq K}\frac{\ep^2[1-\exp(-\frac{2|\bk|^{2\bt}s}{\ep^2})]}{|\bk|^{2\bt}}
\times{|\bk|^2\over
    |\bk|^{2\al+d-d}
\left[(\ep^2+|\bk|^{2\bt})^2+k_1^2\right]}~d\bk~ ds\right]^{1/2}.
\label{1.50}
\ee
The expression (\ref{1.50}) tends to $0$ by virtue of the same change of
variables as the one used above and  the Dominated Convergence Theorem.

Finally, we consider the second term of (\ref{31.12}). After  a simple
 change of variables it can be rewritten as
\be
\label{223}
\frac{2}{\ep}\int\limits_0^tds\int\limits_0^sd{s_1}\int\limits_{s_1}^{s}
\bE\left\{\left(\vV\left(\frac{{s_1}}{\ep^{2}},
\vv\frac{s_1}{\ep^2}+\by_\ep({s_1})\right)\cdot
\nabla\bU_{\ep^2}\left(\frac{{s_2}}{\ep^{2}},\by_\ep({s_1});
\frac{s_1}{\ep^2}\right)\right)
\cdot\bU_{\ep^2}(0,\bze)\right\}ds_2.
\ee

An analogous argument to the one we used to deal with
(\ref{01021}) shows that  (\ref{223}) vanishes as $\ep\da$, provided
that
\be
\label{0330}
{1\over \ep^2}\bE\left|\int\limits_{s_1}^s
\nabla\bU_{\ep^2}\left(\frac{{s_2}}{\ep^{2}},\by_\ep({s_1});
\frac{s_1}{\ep^2}\right)ds_2\right|^2=
{1\over \ep^2}\bE\left|\int\limits_{0}^{s-s_1}
\nabla\bU_{\ep^2}\left(\frac{{s_2}}{\ep^{2}},\bze;0\right)ds_2\right|^2
\ee
\[
\leq{2\over \ep^2}\bE\left|\int\limits_{0}^{s-s_1}
\nabla\bU_{\ep^2}^{(1)}\left(\frac{{s_2}}{\ep^{2}},\bze;0\right)ds_2\right|^2
+{2\over \ep^2}\bE\left|\int\limits_{0}^{s-s_1}
\nabla\bU_{\ep^2}^{(2)}\left(\frac{{s_2}}{\ep^{2}},\bze;0\right)ds_2\right|^2
\]
 tends to zero as $\ep\da$. Here
\be
\label{4163}
\nabla\bU_{\ep^2}^{(i)}(\frac{{s_2}}{\ep^{2}},\bze;0):=\nabla\left[
\vV(\frac{{s_2}}{\ep^{2}},\vv\frac{{s_2}}{\ep^{2}};0)\cdot
\nabla\chi_{\ep^2}^{(i)}(\frac{{s_2}}{\ep^{2}},\bze;0)\right]~~i=1,2.
\ee
Using  spectral representations (\ref{4161}), (\ref{4162}), (\ref{41610}) and
(\ref{981})  we can write  the $p,q$-th entry of the matrix defined by
(\ref{4163}) as being equal to
\[
(-1)^i\int\exp\left\{-\frac{|\bk|^{2\bt}s_2}{\ep^2}\right\}C_i(\bk,\ep^2)
\bk
\hat{V}_{1-i,p}\left(0,\vv\frac{s_2}{\ep^2},d\bk\right)\cdot\int
l_q\exp\left\{-\frac{|\bl|^{2\bt}s_2}{\ep^2}\right\}
\hat{\vV}_{1}\left(0,\vv\frac{s_2}{\ep^2},d\bl\right)
\]
\[
-\int\exp\left\{-\frac{|\bk|^{2\bt}s_2}{\ep^2}\right\}C_i(\bk,\ep^2)\bk
k_q\hat{V}_{i,p}\left(0,\vv\frac{s_2}{\ep^2},d\bk\right)\cdot\int
\exp\left\{-\frac{|\bl|^{2\bt}s_2}{\ep^2}\right\}
\hat{\vV}_{0}\left(0,\vv\frac{s_2}{\ep^2},d\bl\right).
\]
Hence, each term of the rightmost part of (\ref{0330}) equals
\[
\frac{2}{\ep^2}\int\limits_0^{s-s_1}ds_2\int\limits_0^{s-s_1}ds_2'
\int\int\int\int\exp\left\{-\frac{(|\bk|^{2\bt}+|\bl|^{2\bt})s_2}{\ep^2}\right\}
\exp\left\{-\frac{(|\bk'|^{2\bt}+|\bl'|^{2\bt})s_2}{\ep^2}\right\}
\]
\[
\times
C_i(\bk,\ep^2)C_i(\bk',\ep^2) W(d\bk,d\bl,d\bk',d\bl';s_2,s_2')
\]
with $W$ the signed measure given by
\[
\sum\limits_{p,q=1}^d\bE\left\{\left[(-1)^i
k_q\hat{V}_{1-i,p}\left(0,\vv\frac{s_2}{\ep^2},d\bk\right)\bk\cdot
\hat{\vV}_{1}\left(0,\vv\frac{s_2}{\ep^2},d\bl\right)-
k_q\hat{V}_{i,p}\left(0,\vv\frac{s_2}{\ep^2},d\bk\right)\bk\cdot
\hat{\vV}_{0}\left(0,\vv\frac{s_2}{\ep^2},d\bl\right)\right]\right.
\]
\[
\times\left.\left[(-1)^i
k_q'\hat{V}_{1-i,p}\left(0,\vv\frac{s_2'}{\ep^2},d\bk'\right)\bk'\cdot
\hat{\vV}_{1}\left(0,\vv\frac{s_2'}{\ep^2},d\bl'\right)-
k_q'\hat{V}_{i,p}\left(0,\vv\frac{s_2'}{\ep^2},d\bk'\right)\bk'\cdot
\hat{\vV}_{0}\left(0,\vv\frac{s_2'}{\ep^2},d\bl'\right)\right]\right\}.
\]
The above expectation is to be calculated treating the spectral measures as
formal Gaussian random variables and subsequently applying the rules of
calculating the expectation of products of Gaussians with the following
relations for the covariance
\[
\bE\left[\widehat{\vV}_{l}(0,\bx,d\bk)\otimes
\widehat{\vV}_{l'}(0,\bx',d\bk')\right]=
\delta_{l,l'}c_0(\bk\cdot(\bx-\bx'))\hat{\bR}(\bk)
\delta(\bk-\bk')d\bk d\bk'.
\]
After  a straightforward computation we obtain that
\[
W(d\bk,d\bl,d\bk',d\bl';s_2,s_2')=2|\bk|^2c_0(k_1\frac{s_2-s_2'}{\ep^2})
c_0(l_1\frac{s_2-s_2'}{\ep^2})\mbox{tr}\hat{\bR}(\bk)
\bk^{T}\hat{\bR}(\bl)\bk
\]
\[
\times\delta(\bk-\bk')
\delta(\bl-\bl')d\bk d\bl d\bk'd\bl'.
\]

Thus each term on the utmost left hand side of (\ref{0330}) can be bounded from
above by
\be
\label{034}
C\mathop{\int\int}\limits_{|\bk|,|\bl|\leq K}
C_i^2(\bk,\ep^2)\frac{|\bk|^4\ep^2\left[1-\exp\left\{-{(|\bk|^{2\bt}+|\bl|^{2\bt})(s-s_1)\over
    \ep^2}\right\}\right]^2}{(|\bk|^{2\bt}+|\bl|^{2\bt})^2(|\bk||\bl|)^{2\al+d-2}}
~d\bk~ d\bl
\ee
with the constant $C$ independent of $\ep$. For $\al+2\bt<2$
 (\ref{034})
tends to zero, as $\ep\da$, by virtue of the Dominated Convergence Theorem and
the proof of the Lemma \ref{lemma1} is finished\kropa

{\em The proof of Theorem \ref{theorem1}.} By virtue of (\ref{124}) we write
\be
\label{90.10}
\by_\ep(t)=M_\ep(t)+R_\ep(t)
\ee
where the remainder term $R_\ep(t)$ consists of a stationary, i.e.
$-\chi_\ep(t)$, and an additive part
\[
R_{a,\ep}(t):=\int\limits_0^t\chi_{\ep}(s)ds+
\int\limits_0^t\bU_\ep(s)ds.
\]

According to Lemma \ref{lemma1}, $M_\ep(t)$, $t\geq0$, is a Brownian motion
 whose covariance
 matrix, given by (\ref{061104}), converges to $\bD_0$ as
 $\ep\da$. In consequence, $M_\ep(t)$, $t\geq0$, converges weakly, when $\ep\da$,
   as processes with continuous trajectories (cf. e.g. \cite{helland}).

In order to conclude the proof of Theorem \ref{theorem1} it remains to be seen
 that for
 arbitrary $T>0$, $\delta>0$ we have
\[
\lim\limits_{\ep\da}P[\sup\limits_{0\leq t\leq T}|R_\ep(t)|\geq
\delta]=0.
\]
The assertions of Lemma \ref{lemma1} imply that we  only need to show
\[
\lim\limits_{\ep\da}P[\sup\limits_{0\leq t\leq T}|\chi_{\ep}(t)|\geq
\delta]=0.
\]
Since $\chi_{\ep}$ is stationary we can write
\be
\label{061107}
P[\sup\limits_{0\leq t\leq T}|\chi_{\ep}(t)|\geq
\delta]\leq\frac{T}{\ep^2}P[\sup\limits_{0\leq t\leq T\ep^2}|\chi_{\ep}(t)|\geq
\delta].
\ee
Using again (\ref{90.10}) to represent $\chi_{\ep}$ we can estimate
the right hand side of (\ref{061107}) by the sum
\be
\label{061110}
\frac{T}{\ep^2}P[\sup\limits_{0\leq t\leq T\ep^2}
[|R_{a,\ep}(t)|+|\by_\ep(t)|]\geq
\delta]+\frac{T}{\ep^2}P[\sup\limits_{0\leq t\leq T\ep^2}
|M_\ep(t)|\geq
\delta]+
\frac{T}{\ep^2}P[|\chi_{\ep}(0)|\geq
\delta].
\ee
The last term in (\ref{061110}) can be estimated using Chebychev's inequality by
\be
\label{061112}
\frac{\bE|\chi_{\ep}(0)|^pT}{\ep^2\delta^p}.
\ee
Choosing $p$ so that $\ga\cdot  p>4$ with
 $\ga$  given by Lemma \ref{lemma1}, we see  that the expression
(\ref{061112}) tends to $0$ as
$\ep\da$ by virtue of  (\ref{08121}).

The second term of  (\ref{061110}) vanishes as
$\ep\da$ thanks to
\[
\lim\limits_{\ep\da}\ep^{-2}\bE|M_\ep(\ep^2T)|^2=0
\]
and
 an elementary martingale inequality
\be
\label{08142}
\frac{T}{\ep^2}P[\sup\limits_{0\leq t\leq T\ep^2}|M_\ep(t)|]\geq
\delta]\leq\frac{\bE|M_\ep(\ep^2T)|^2T}{\ep^2\delta^2}.
\ee

The first term of (\ref{061110})  can be estimated by
\be
\label{061111}
C\left\{\frac{\ep^2T}{\delta^2}\bE|\chi_{\ep}(0)|^2+
\frac{\ep^2T}{\delta^2}\bE|\bU_{\ep}(0)|^2+
\frac{T}{\ep^2}P\left[\int\limits_0^T
\left|\vV(s,\vv s+\by(s))\right|ds\geq\frac{\delta}{\ep}\right]\right\}.
\ee
The  first two expressions in (\ref{061111}) tend to $0$
thanks to Lemma
\ref{lemma1}. The last
term of (\ref{061111})
vanishes as $\ep\da$, provided that $p>2$. This is due to the fact that
\[
\limsup\limits_{\ep\da}\bE\left[\int\limits_0^T
\left|\vV(s,\vv s+\by(s))\right|ds\right]^p\leq
\bE\left|\vV(0,\bze)\right|^pT^{p}<+\infty.
\]
This finishes the proof of Theorem \ref{theorem1}~\kropa

\section{Proof of Theorem \ref{theorem2}.} \label{section2}
Throughout this section we  assume that the scaled
 trajectory $\by_\ep$ is given by (\ref{120801})
with $\delta_{\al,\bt}$ given by (\ref{01051}), (\ref{01052}). To avoid
cumbersome notation we  suppress writing the subscript of $\delta$.

Define the path corrector and convector by
\[
\chi_\ep(t):=\ep\chi_{\ep^{2\delta}}(\frac{t}{\ep^{2\delta}},
\by_\ep(t)),
\]
and
\[
\bU_\ep(t):=\bU_{\ep^{2\delta}}(\frac{t}{\ep^{2\delta}},\by_\ep(t)).
\]
The following lemma can be proven exactly the same way as Lemma \ref{lemma1}.
\begin{lemma}
\label{lemma61}
The scaled trajectory of a particle satisfies
\be
\label{7.10}
\by_\ep(t)=-\chi_\ep(t)+\chi_\ep(0)+\int\limits_0^t\chi_\ep(s)ds+
\ep^{2(1-\delta)}\int\limits_0^t\bU_{\ep}(s)ds+\ep^{1-\delta}M_{\ep}(t),
\ee
where $M_{\ep}(t)$ is exactly as in the statement of Lemma {\em \ref{lemma1}}.
\end{lemma}

{\bf Remark.} Calculations analogous to those made after
 (\ref{281101}) show that for $\bt\leq1/2$
the exponent $\delta$ is chosen to  make $\bE|\ep^{1-\delta}M_{\ep}(t)|^2$ of
order $O(1)$. Indeed,
\be
\label{01053}
\bE|\ep^{1-\delta}M_{\ep}(t)|^2=
\ep^{2-2\delta} \mbox{tr}\bD_\ep t\simeq
\ep^{2\left(1+\delta\frac{1-\al-2\bt}{\bt}\right)}
\ee
and
\be
\label{01054}
\bE|\chi_\ep(t)|^2=\ep^2\bE|\chi_{\ep^{2\delta}}(0,\bze)|^2
\simeq\ep^{2\left(1+\delta\frac{1-\al-2\bt}{\bt}\right)}.
\ee
In the case when $\bt>1/2$ we can deduce by the same argument that the terms on
the left hand sides of (\ref{01053}) and (\ref{01054}) vanish as $\ep\da$.

We note here that $\delta<1$ for $\alpha+\beta>1$ and $ \beta>0$. Thanks to the
assumption that $\al+2\bt<2$ we have (cf. (\ref{40'}))
\be
\label{08171}
\ep^{4-4\delta}\bE|\nabla \chi_{\ep^{2\delta}}(0,\bze)|^2\sim
\ep^{4-4\delta}\da\qquad\mbox{as}~~\ep\da.
\ee
This together with (\ref{221}) implies that the term on the right hand side of
(\ref{7.10}) involving the convector vanishes in probability as $\ep\da$. In
addition, we deduce by (\ref{08171}) that
\be
\label{061705}
\chi_\ep(t)=\ep\chi_{\ep^{2\delta}}\left(\frac{t}{\ep^{2\delta}},\bze\right)+
\ep^{2(1-\delta)}\int\limits_0^t\nabla\chi_{\ep^{2\delta}}
\left(\frac{s}{\ep^{2\delta}},\by_\ep(s)\right)\cdot
\vV\left(\frac{s}{\ep^{2\delta}},\vv\frac{s}{\ep^{2\delta}}+
\by_\ep(s)\right)ds
\ee
\[
=\ep\chi_{\ep^{2\delta}}\left(\frac{t}{\ep^{2\delta}},\bze\right)+o_\ep(t).
\]
Here and in the sequel we use the generic notation $o_\ep(t)$ for any stochastic
process which satisfies $\lim\limits_{\ep\da}\bE\sup\limits_{0\leq t\leq
T}|o_\ep(t)|^2=0$, for any $T>0$.

With this, (\ref{7.10}) becomes
\be
\label{18.10}
\by_\ep(t)=-\ep\chi_{\ep^{2\delta}}\left(\frac{t}{\ep^{2\delta}},\bze\right)+
\ep\chi_{\ep^{2\delta}}(0,\bze)+\ep\int\limits_0^t
\chi_{\ep^{2\delta}}\left(\frac{s}{\ep^{2\delta}},\bze\right)ds+
\ep^{1-\delta}M_{\ep}(t)+o_\ep(t).
\ee

\begin{lemma}
\label{lemma62}
We have the following expansion of the corrector:
\be
\label{061715}
\ep\chi_{\ep^{2\delta}}\left(\frac{t}{\ep^{2\delta}},\bze\right)=
\ep\chi_{\ep^{2\delta}}(0,\bze)+\ep\int
\limits_0^{t/\ep^{2\delta}}\chi_{\ep^{2\delta}}(s,\bze)ds-
\ep\int\limits_0^{t/\ep^{2\delta}}\vV(s,\vv s)ds+
\ep^{1-\delta}\tilde{M}_{\ep}(t)
\ee
where $\tilde{M}_{\ep}$ is a Brownian motion which satisfies
\begin{itemize}
\item[{\em i)}] the $i,j$-th entry of its covariance
 matrix is given by {\em (\ref{a1.10})}.
\item[{\em ii)}] there exists $\ga>0$ such that
\be
\label{061711}
\lim\limits_{\ep\da}\ep^{2(1-\delta)-\ga}
\bE\left|M_{\ep}(t)-\tilde{M}_{\ep}(t)\right|^2=0.
\ee
\end{itemize}
\end{lemma}

{\bf Proof.} The proofs of (\ref{061715}) and part i) go along the lines of the
proofs of the corresponding parts of Lemma \ref{lemma1} so we leave this
argument out. Actually, the fact that the spatial argument is fixed
significantly simplifies
 the proof in the present case.

We focus on the proof of part ii). We first calculate the pseudogenerator
\[
\cL(t):=\ep^{2(1-\delta)}\cL\left|M_{\ep}(t)-\tilde{M}_{\ep}(t)\right|^2.
\]
In light of (\ref{061711}) and (\ref{7.10}) it is clear that
\be
\label{18.11}
\cL(t)=\ep^2\cL|\bZ(t,\by_\ep(t))|^2-\ep^2 \bZ(t,\by_\ep(t))\cdot
\cL \bZ(t,\by_\ep(t))
\ee
where
\[
\bZ(t,\by):=\chi_{\ep^{2\delta}}\left(\frac{t}{\ep^{2\delta}},\by\right)-
\chi_{\ep^{2\delta}}\left(\frac{t}{\ep^{2\delta}},\bze\right).
\]
The first term on the right hand side of (\ref{18.11}) can
be computed as $L^1$ limit of
\[
\ep^2\lim\limits_{r\da}\frac{\bE\left\{|\bZ(t+r,\by_\ep(t))|^2
\left|\right.{\cal
      V}_{-\infty,t/\ep^{2\delta}}\right\}-
|\bZ(t,\by_\ep(t))|^2}{r}.
\]
A direct computation similar to that performed after (\ref{010516})
shows that
\be
\label{010510}
\cL(t)=2\ep^{2(1-\delta)}\int\limits_{R^d}
\frac{|\bk|^{2\bt}a(|\bk|)}{|\bk|^{2\al+d-2}
[(\ep^{2\delta}+|\bk|^{2\bt})^2+k_1^2]} [1-c_0(\bk\cdot\bx_\ep(t))]~d\bk.
\ee

Here, recall, $\bx_\ep(t)=\vv t+\by_\ep(t)$. For $2\bt>1$ computations analogous
to those performed after (\ref{281101}) or (\ref{40'}) convince us that the
right hand side of (\ref{010510}) vanishes as $\ep\da$.

To estimate the mathematical expectation of the right hand side of
 (\ref{010510}) when $2\bt\leq1$ we break it
into two disjoint integrals: the first over $|\bk|\leq
\ep^{\delta\varrho/\bt}$, with $0<\varrho<1$, and the second over its
complement.

The expectation of the second integral can be
estimated by
\be
\label{010511}
2\ep^{2(1-\delta)}\int\limits_{K\geq|\bk|\geq
\ep^{\delta\varrho/\bt}}\frac{|\bk|^{2\bt+2-2\al-d}}
{[(\ep^{2\delta}+|\bk|^{2\bt})^2+k_1^2]}~d\bk.
\ee
Employing again the procedure used to estimate the last part in (\ref{281101})
we infer that (\ref{010511}) is less than or equal to a constant times
\[
\ep^{2(1-\delta)}\int\limits_{\ep^{\frac{\delta\varrho}{\bt}}}^K
\frac{r^{2\bt+1}}{r^{2\al}(r^{2\bt}+\ep^{2\delta})^2}~dr.
\]
With the change of variables $r':=\ep^{{\delta}{\bt}}r$ we conclude that the
above expression vanishes, as $\ep\da$, at the rate $O(\ep^\ga)$ for some
$\ga>0$ since $0<\varrho<1$.

The expectation of the first  integral can be estimated as follows:
\beq
&&
\label{01.21}
2\ep^{3-4\delta}\int \limits_0^t ds\int\limits_{|\bk|\leq
  \ep^{\delta\varrho/\bt}}\frac{|\bk|^{2\bt}
a(|\bk|)}{|\bk|^{2\al+d-2}[(\ep^{2\delta}+|\bk|^{2\bt})^2+k_1^2]}
\bE\left[-\sin(\bk\cdot\bx_\ep(s))\bk\cdot\vV\left(\frac{s}{\ep^{2\delta}},
\bx_\ep(s)\right)\right]~d\bk
\nonumber\\
&&\leq C\ep^{3-4\delta}\int\limits_{|\bk|\leq\ep^
{\delta\varrho/\bt}}\frac{|\bk|^{2\bt+3-2\al-d}}
{[(\ep^{2\delta}+|\bk|^{2\bt})^2+k_1^2]}d\bk.
\eeq

Via the same type of argument as used to estimate the previous integral we find
that (\ref{01.21}) is bounded from above by
\[
\ep^{\frac{2}{\al+2\bt-1}}\int
\limits_0^{\ep^{\delta(\varrho-1)/\bt}}
\frac{r^{2\bt+2}}{r^{2\al}(r^{2\bt}+1)^2}dr.
\]
The divergence of the integral in (\ref{01.21}), if occurs, can be made
arbitrarily slow by choosing  $\varrho$ sufficiently close to one. Thus  the
expression in (\ref{01.21}) vanishes for $\alpha+2\beta\leq2$. This concludes
the proof of Lemma \ref{lemma62}\kropa

As a consequence of Lemma \ref{lemma62} we obtain
\be
\label{08172}
\bx_\ep(t)=\by_\ep(t)+o_\ep(t),
\ee
where
\be
\label{062502}
\by_\ep(t)=\ep\int^{t/\ep^{2\delta}}_0\vV(s,\bze)ds
\ee
is a Gaussian process whose covariance matrix equals
\be
\label{010512}
\bE\left[\by_\ep(t)\otimes\by_\ep(t)\right]=2\ep^2\int\limits_0^
{t/\ep^{2\delta}}ds\int\limits_0^{s}\bR(s',\vv s')~ds'.
\ee
An elementary calculation  shows that the right hand side of (\ref{010512})
tends to $\bD t$ as $\ep\da$ with $\bD$ given by (\ref{010513}). Hence $\by_\ep$
converges weakly as $\ep\da$ to a fractional Brownian motion $\bB_H(t)$ whose
covariance is given by (\ref{010514}). \\
{\small Authors Addresses:\\ Albert
Fannjiang, Department of Mathemathics, University of California,  Davis, CA
95616-8633, USA\\Tomasz Komorowski, Instytut Matematyki, UMCS, pl. Marii-Curie
Sk\l odowskiej 1, 20-031 Lublin, Polska}

\end{document}